\documentclass[11pt]{amsart}
\usepackage{palatino}
\usepackage{amsfonts}
\usepackage{amssymb}
\usepackage{amscd}
\usepackage[breaklinks,bookmarksopen,bookmarksnumbered]{hyperref}
\usepackage[dvips]{graphicx}

\newtheorem{theorem}{Theorem}[section]
\newtheorem{proposition}[theorem]{Proposition}
\newtheorem{corollary}[theorem]{Corollary}

\theoremstyle{definition}
\newtheorem{definition}[theorem]{Definition}
\newtheorem{remark}[theorem]{Remark}

\newtheorem{conjecture/question}[theorem]{Conjecture/Question}

\newtheorem{remark/definition}[theorem]{Remark/Definition}
\newtheorem{terminology/notation}[theorem]{Terminology/Notation}

\def\PP{{\textbf P}}
\def\OO{\mathcal{O}}

\def\cB{\mathcal{B}}
\def\cA{\mathcal{A}}

\def\P{\mathcal{P}}
\def\W{\mathcal{W}}

\def\cS{\mathcal{S}}

\def\L{\mathcal{L}}

\def\cM{\mathcal{M}}

\def\rr{\overline{\mathcal{R}}}
\def\cZ{\mathcal{Z}}
\def\cU{\mathcal{U}}
\def\cC{\mathcal{C}}
\def\H{\mathcal{H}}
\def\Pic0{{\rm Pic}^0(X)}

\def\mm{\overline{\mathcal{M}}}
\def\ss{\overline{\mathcal{S}}}

\def\dd{\overline{\mathcal{D}}}

\def\aux{\mathfrak G^1_k(\overline{\textbf{S}}_g^0/\rem_g^0)}
\def\thet{\overline{\Theta}_{\mathrm{null}}}

\def\AUX{\mathfrak G^r_d(\overline{\textbf{S}}_g^0/\rem_g^0)}
\def\ts{\overline{\textbf{S}}}
\def\pem{\widetilde{\textbf{M}}}
\def\pes{\widetilde{\textbf{S}}}
\def\rem{\overline{\textbf{M}}}

\def\tet{\overline{\Theta}}

\begin{document}
\title{Brill-Noether geometry on moduli spaces of spin curves}

\author[G. Farkas]{Gavril Farkas}

\address{Humboldt-Universit\"at zu Berlin, Institut F\"ur Mathematik,  Unter den Linden 6
\hfill \newline\texttt{}
 \indent 10099 Berlin, Germany} \email{{\tt farkas@math.hu-berlin.de}}
\thanks{}

\thanks{Research partially supported by Sonderforschungsbereich  "Raum-Zeit-Materie".}

\maketitle

The aim of this paper is to initiate a study of geometric divisors of Brill-Noether type on the moduli space $\ss_g$ of spin curves of genus $g$. The moduli space $\ss_g$ is a compactification the parameter space $\cS_g$ of  pairs $[C, \eta]$, consisting of a smooth genus $g$ curve $C$ and a theta-characteristic $\eta\in \mbox{Pic}^{g-1}(C)$, see \cite{C}. The study of the birational properties of $\ss_g$ as well as other moduli spaces of curves with level structure has received an impetus in recent years, see \cite{BV} \cite{FL}, \cite{F2}, \cite{Lud}, to mention only a few results. Using syzygy divisors, it has been proved in \cite{FL} that the Prym moduli space $\rr_g:=\mm_g(\mathcal{B} \mathbb Z_2)$ classifying curves of genus $g$ together with a point of order $2$ in the Jacobian variety, is a variety of general type for $g\geq 13$ and $g\neq 15$.
The moduli space $\ss_g^+$ of even spin curves of genus $g$ is known to be of general type for $g>8$, uniruled for $g<8$, see \cite{F2}, whereas the Kodaira dimension of $\ss_8^+$ is equal to zero, \cite{FV}. This was the first example of a naturally defined moduli space of curves of genus $g\geq 2$, having intermediate Kodaira dimension. An application of the main construction of this paper, gives a new way of computing the class of the divisor $\thet$ of vanishing theta-nulls on $\ss_g^+$, reproving thus the main result of \cite{F2}.

Virtually  all attempts to show that a certain moduli space $\mm_{g, n}$ is of general type, rely on the calculation of certain effective divisors $D\subset \mm_{g, n}$ enjoying extremality properties in their effective cones $\mbox{Eff}(\mm_{g, n})$, so that the canonical class $K_{\mm_{g, n}}$ lies in the cone spanned by $[D]$, boundary classes
$\delta_{i: S}$, tautological classes $\lambda, \psi_1, \ldots, \psi_n$, and possible other effective geometric classes. Examples of such a program being carried out, can be found in \cite{EH2}, \cite{HM}-for the case of \emph{Brill-Noether divisors} on $\mm_g$ consisting of curves  with a $\mathfrak g^r_d$ when $\rho(g, r, d)=-1$, \cite{Log}-where pointed Brill-Noether divisors on $\mm_{g, n}$ are studied, and \cite{F1}-for the case of \emph{Koszul divisors} on $\mm_g$, which provide counterexamples to the Slope Conjecture on $\mm_g$. A natural question is what the analogous geometric divisors on the spin moduli space of curves $\ss_g$ should be?

\vskip 6pt
In this paper we propose a construction for \emph{spin Brill-Noether divisors} on both spaces $\ss_g^+$ and $\ss_g^-$, defined in terms of the relative position of theta-characteristics with respect to difference varieties on Jacobians. Precisely, we fix integers $r, s\geq 1$ such that $d:=rs+r \equiv 0 \mbox{ mod }2$, and then set $g:=rs+s$. One can write $d=2i$. By standard Brill-Noether theory, a general curve $[C]\in \cM_g$ carries a finite number of (necessarily complete and base point free) linear series $\mathfrak g^r_d$. One considers the following loci of spin curves (both odd and even)
$$\cU_{g, d}^r:=\{[C, \eta]\in \cS_g^{\mp}: \exists L\in W^r_d(C)\ \mbox{ such that } \ \eta\otimes L^{\vee}\in C_{g-i-1}-C_i\}.$$
Thus $\cU_{g, d}^r$ consists of spin curves such that the embedded curve $C\stackrel{|\eta\otimes L|}\longrightarrow \PP^{d-1}$ admits an $i$-secant $(i-2)$-plane.
We shall prove that for $s\geq 2$, the locus $\cU_{g, d}^r$ is always a divisor on $\cS_{g}^\mp$, and we find a formula for the class of its compactification in $\ss_g^{\mp}$. For simplicity, we display this formula in the introduction only in the case $r=1$, when $g\equiv 2 \mbox{ mod }4$:
\begin{theorem}\label{spinclass}
We fix an integer $a\geq 1$ and set $g:=4a+2$. The locus
$$\cU_{4a+2, 2a+2}^{1}:=\{[C, \eta]\in \cS_{4a+2}^{\mp}: \exists L\in W_{2a+2}^{1}(C) \mbox{ such that } \eta\otimes L^{\vee} \in C_{3a}-C_{a+1}\}$$
is an effective divisor and the class of its compactification in $\ss_{g}^{\mp}$ is given by
$$\overline{\cU}_{4a+2, 2a+2}^{1}\equiv {4a\choose a}{4a+2\choose 2a}\frac{a+2}{8(2a+1)(4a+1)}\Bigl(\bigl(192a^3+736a^2+692a+184\bigr)\lambda-$$
$$-\bigl(32a^3+104a^2+82a+19\bigr)\alpha_0-\bigl(64a^3+176a^2+148a+36\bigr)\beta_0-\cdots \Bigr)\in \mathrm{Pic}(\ss_g^{\mp}).$$
\end{theorem}
To specialize further, in Theorem \ref{spinclass} we set $a=1$, and find the class of (the closure of) the locus of spin curves $[C, \eta]\in \cS_6^{\mp}$, such that
there exists a pencil $L\in W^1_4(C)$ for which the linear series $C\stackrel{|\eta\otimes L|}\longrightarrow \PP^3$ is not very ample:
$$\overline{\cU}_{6, 4}^1\equiv 451 \lambda-\frac{237}{4}\alpha_0-106\beta_0-\cdots\in \mathrm{Pic}(\ss_6^{\mp}\bigr).$$

The case $s=1$, when necessarily $L=K_C\in W^{g-1}_{2g-2}(C)$, produces a divisor only on $\ss_g^+$, and we recover in this way the main calculation
from \cite{F2}, used to prove that $\ss_g^+$ is a variety of general type for $g>8$. We recall that $\Theta_{\mathrm{null}}:=\{[C, \eta]\in \cS_g^+: H^0(C, \eta)\neq 0\}$ denotes the divisor of \emph{vanishing theta-nulls}.
\begin{theorem}\label{f2}
Let $\pi:\ss_g^+\rightarrow \mm_g$ be the ramified covering which forgets the spin structure. For $g\geq 3$, one has the following equality $\overline{\mathcal{U}}_{g, 2g-2}^{g-1}=2\cdot \thet$ of codimension $1$-cycles on the open subvariety $\pi^{-1}(\cM_g\cup \Delta_0)$ of $\ss_g^+$. Moreover, there is an equality of classes
$$\overline{\mathcal{U}}_{g, 2g-2}^{g-1}\equiv 2\cdot \overline{\Theta}_{\mathrm{null}}\equiv \frac{1}{2}\lambda-\frac{1}{8}\alpha_0-0\cdot \beta_0-\cdots \ \in \mathrm{Pic}(\ss_g^+).$$
\end{theorem}
We remark once more, the low slope of the divisor $\thet$. No similar divisor with such remarkable class is known to exist on $\rr_g$.  In Section 4, we present a third way of calculating the class $[\thet]$, by rephrasing the condition that a curve $C$ have a vanishing theta-null $\eta$, if and only if, for a pencil $A$ on $C$ of minimal degree, the multiplication map of sections
$$H^0(C, A)\otimes H^0(C, A\otimes \eta)\rightarrow H^0(C, A^{\otimes 2}\otimes \eta)$$
is not an isomorphism. For $[C]\in \cM_g$ sufficiently general, we note that $$\mbox{dim } H^0(C, A)\otimes H^0(C, A\otimes \eta)=\mbox{dim } H^0(C, A^{\otimes 2}\otimes \eta).$$ In this way, $\thet$ appears as the push-forward of a degeneracy locus of a morphism between vector bundles of the same rank defined over a Hurwitz stack of coverings. To compute the push-forward of tautological classes from a Hurwitz stack, we use the techniques developed in \cite{F1} and \cite{Kh}.

In the last section of the paper, we study the divisor $\tet_{g, 1}$ on the universal curve $\mm_{g, 1}$, which consists of points in the support of odd theta-characteristics. This divisor, somewhat similar to the divisor $\overline{\W}_g$ of Weierstrass points on $\mm_{g, 1}$, cf. \cite{Cu}, should be of some importance in the study of the birational geometry of $\mm_{g, 1}$:

\begin{theorem}\label{oddthet} The class of the compactification in $\mm_{g, 1}$ of the effective divisor
$$\Theta_{g, 1}:=\{[C, q]\in \cM_{g, 1}: q\in \mathrm{supp}(\eta)\ \mbox{ for some }\ [C, \eta]\in \cS_g^-\}$$ is given by the following formula:
$$\tet_{g, 1}\equiv 2^{g-3}\Bigl((2^g-1)\bigl(\lambda+2\psi\bigr)-2^{g-3}\delta_{\mathrm{irr}}-(2^g-2)\delta_1- \sum_{i=1}^{g-1}(2^i+1)(2^{g-i}-1)\delta_i\Bigr) \in \mathrm{Pic}(\mm_{g, 1}).$$
\end{theorem}
When $g=2$, the divisor $\Theta_2$ specializes to the divisor of Weierstrass points:
$$\Theta_{2, 1}=\W_2:=\{[C, q]\in \cM_{2, 1}: q\in C \ \mbox{ is a Weierstrass point}\}.$$
If we use Mumford's formula $\lambda=\delta_0/10+\delta_1/5\in \mathrm{Pic}(\mm_2)$, Theorem \ref{oddthet} reads
$$\tet_{2, 1}\equiv \frac{3}{2}\lambda+3\psi-\frac{1}{4}\delta_{\mathrm{irr}}-\frac{3}{2}\delta_1=-\lambda+3\psi-\delta_1\in \mathrm{Pic}(\mm_{2, 1}),$$
that is, we recover the formula for the class of the Weierstrass divisor on $\mm_{2, 1}$, cf. \cite{EH2}. When $g=3$,
 the condition $[C, q]\in \Theta_{3, 1}$, states that the point $q\in C$ lies on one of the  $28$ bitangent lines of the canonically embedded curve
$C\stackrel{|K_C|}\longrightarrow \PP^2$.
\begin{corollary}
The class of the compactification in $\mm_{3, 1}$ of the bitangent locus $$\Theta_{3, 1}:=\{[C, q]\in \cM_{3, 1}: q\  \mbox{ lies on a bitangent of } C\}$$
is equal to \
$\tet_{3, 1} \equiv 7\lambda+14\psi-\delta_{\mathrm{irr}}-9\delta_1-5\delta_2\in \mathrm{Pic}(\mm_{3, 1})$.
\end{corollary}
If $p:\mm_{g, 1}\rightarrow \mm_g$ is the map forgetting the marked point, we note the equality
$$\dd_3\equiv p^*(\mm_{3, 2}^1)+2\cdot\overline{\W}_3+2\psi \in \mbox{Pic}(\mm_{3, 1}),$$ where $\overline{\W}_3\equiv -\lambda+6\psi-3\delta_1-\delta_2$ \ is the divisor of Weierstrass points on $\mm_{3, 1}$.
Since the class $\psi\in \mbox{Pic}(\mm_{3, 1})$ is big and nef, it follows that $\tet_{3, 1}$ (unlike the divisor $\tet_{2, 1} \in \mbox{Pic}(\mm_{2, 1})$), lies in the interior of the cone of effective divisors $\mbox{Eff}(\mm_{3, 1})$, or it other words, it is big. In particular, it cannot be contracted by a rational map $\mm_{3, 1}\dashrightarrow X$ to any projective variety $X$. This phenomenon extends to all higher genera:
\begin{corollary}
For every $g\geq 3$, the divisor $\tet_{g, 1}\in \mathrm{Eff}(\mm_{g, 1})$ is big.
\end{corollary}
It is not known whether the Weierstrass divisor $\overline{\mathcal{W}}_g$ lies on the boundary of the effective cone $\mbox{Eff}(\mm_{g, 1})$ for $g$ sufficiently large.

\section{Generalities about $\ss_g$}

As usual, we follow that the convention that if $\textbf{M}$ is a Deligne-Mumford stack, then $\cM$ denotes its associated coarse moduli space.  We first recall basic facts about Cornalba's stack of stable spin curves
$\pi:\ts_g\rightarrow \rem_g$, see \cite{C}, \cite{F2}, \cite{Lud} for details and other basic properties. If $X$ is a nodal curve,
a smooth rational component $R\subset X$ is said to be
\emph{exceptional} if $\#(R\cap \overline{X-R})=2$. The curve $X$ is
said to be \emph{quasi-stable} if $\#(R\cap \overline{X-R})\geq 2$ for any smooth rational component $R\subset X$, and moreover, any two exceptional components of $X$ are disjoint. A quasi-stable curve is obtained from a stable
curve by possibly inserting a rational curve at each of its nodes. We denote by $[\mathrm{st}(X)]\in
\mm_g$ the stable model of the quasi-stable curve $X$.

\begin{definition}\label{spinstructures} A \emph{spin curve} of genus
$g$ consists of a triple $(X, \eta, \beta)$, where $X$ is a genus
$g$ quasi-stable curve, $\eta\in \mathrm{Pic}^{g-1}(X)$ is a line
bundle of degree $g-1$ such that $\eta_{R}=\OO_R(1)$ for every
exceptional component $R\subset X$, and $\beta:\eta^{\otimes
2}\rightarrow \omega_X$ is a sheaf homomorphism which is generically
non-zero along each non-exceptional component of $X$.
\end{definition}

Stable spin curves of genus $g$ form a smooth Deligne-Mumford stack $\ts_g$ which splits into two connected components
$\ts_g^+$ and $\ts_g^-$, according to the parity of $h^0(X, \eta)$. Let $f:\mathcal{C} \rightarrow \ts_g$ be the universal family of spin curves of genus $g$. In particular,  for every point $[X, \eta, \beta]\in \ss_g$, there is an isomorphism between $f^{-1}([X,
\eta, \beta])$ and the quasi-stable curve $X$.
There exists a (universal) spin line bundle $\P \in
\mathrm{Pic}(\mathcal{C})$ of relative degree $g-1$, as well as a morphism of $\OO_{\mathcal C}$-modules
$B:\P^{\otimes 2}\rightarrow \omega_f$ having the property
that $\P_{| f^{-1}([X, \eta, \beta])}=\eta$ and $B_{| f^{-1}([X,
\eta, \beta])}=\beta:\eta^{\otimes 2}\rightarrow \omega_X$, for all
spin curves $[X, \eta, \beta]\in \ss_g$. Throughout we use the canonical isomorphism
$\mbox{Pic}(\bf{\overline{S}}_g)_{\mathbb Q} \cong \mbox{Pic}(\ss_g)_{\mathbb Q}$ and we make little distinction between line bundles on the stack and the corresponding moduli space.

\subsection{The boundary divisors of $\ss_g$}\hfill

We discuss the structure of the boundary divisors of $\ss_g$ and concentrate on the case of $\ss_g^+$, the differences compared to the situation on $\ss_g^-$ being minor. We describe the pull-backs of the boundary divisors $\Delta_i\subset \mm_g$ under the map $\pi$. First we fix an integer $1\leq i\leq [g/2]$ and let $[X, \eta, \beta]\in \pi^{-1}([C\cup_y D])$,
where $[C, y]\in \cM_{i, 1}$ and $[D, y]\in \cM_{g-i, 1}$. For degree reasons, then
$X=C\cup_{y_1} R\cup_{y_2} D$, where $R$ is an exceptional
component such that $C\cap R=\{y_1\}$ and $D\cap R=\{y_2\}$.
Furthermore $\eta=\bigl(\eta_C, \eta_D, \eta_R=\OO_R(1)\bigr)\in
\mbox{Pic}^{g-1}(X)$, where $\eta_C^{\otimes 2}=K_C$ and $\eta_D^{\otimes
2}=K_D$. The theta-characteristics $\eta_C$ and $\eta_D$ have
the same parity in the case of $\ss_g^+$ (and opposite parities for $\ss_g^-$). One denotes by $A_i\subset \ss_g^+$ the closure of
the locus corresponding to pairs of pointed spin curves
$$\bigl([C, y, \eta_C], [D, y, \eta_D]\bigr)\in
\cS_{i, 1}^+\times \cS_{g-i, 1}^+$$ and by $B_i\subset \ss_g^+$ the
closure of the locus corresponding to pairs $$\bigl([C, y, \eta_C], [D,
y, \eta_D]\bigr)\in \cS_{i, 1}^-\times \cS_{g-i, 1}^{-}.$$ If $\alpha:=[A_i], \beta_i:=[B_i]\in \mbox{Pic}(\ts_g^+)$, we have the relation
$\pi^*(\delta_i)=\alpha_i+\beta_i$.

Next, we describe $\pi^*(\delta_0)$ and pick a stable spin curve $[X, \eta, \beta]$ such that $\mathrm{st}(X)=C_{yq}:=C/y\sim q$,
with $[C, y, q]\in \cM_{g-1, 2}$. There are two possibilities
depending on whether $X$ possesses an exceptional component or not.
If $X=C_{yq}$ and $\eta_C:=\nu^*(\eta)$ where $\nu:C\rightarrow X$
denotes the normalization map, then $\eta_C^{\otimes 2}=K_C(y+q)$.
For each choice of $\eta_C\in \mathrm{Pic}^{g-1}(C)$ as above, there
is precisely one choice of gluing the fibres $\eta_C(y)$ and
$\eta_C(q)$ such that $h^0(X, \eta)\equiv 0 \mbox{ mod } 2$. We denote by $A_0$ the
closure in $\ss_g^+$ of the locus of points $[C_{yq}, \ \eta_C\in \mathrm{Pic}^{g-1}(C), \ \eta_C^{\otimes 2}=K_C(y+q)]$ as above.

If $X=C\cup_{\{y, q\}} R$, where $R$ is an exceptional component,
then $\eta_C:=\eta\otimes \OO_C$ is a theta-characteristic on $C$.
Since $H^0(X, \omega)\cong H^0(C, \omega_C)$, it follows that $[C,
\eta_C]\in \cS_{g-1}^{+}$.  We denote by
$B_0\subset \ss_g^+$ the closure of the locus of points
$$\bigl[C\cup_{\{y, q\}} R, \ \eta_C\in \sqrt{K_C}, \ \eta_R=\OO_R(1)\bigr]\in \ss_g^+.$$ A local analysis carried out in \cite{C}, shows that $B_0$ is the branch locus of $\pi$ and the ramification is simple. If
$\alpha_0=[A_0]\in \mbox{Pic}(\ss_g^+)$ and $\beta_0=[B_0]\in
\mbox{Pic}(\ss_g^+)$,  we have the relation
\begin{equation}\label{del0}
\pi^*(\delta_0)=\alpha_0+2\beta_0.
\end{equation}

\section{Difference varieties and theta-characteristics}

We describe a way of calculating the class of a series of effective divisors on both moduli spaces $\ss_g^-$ and $\ss_g^+$, defined in terms of the relative position of a theta-characteristic with respect to the divisorial difference varieties in the Jacobian of a curve. These loci, which should be thought of as divisors of Brill-Noether type on $\ss_g$, inherit a determinantal description over the entire moduli stack of spin curves, via the interpretation of difference varieties in $\mbox{Pic}^{g-2i-1}(C)$ as Raynaud theta-divisors for exterior powers of Lazarsfeld bundles provided in \cite{FMP}. The determinantal description is then extended over a partial compactification $\pes_g$ of $\bf{S}_g$, using the explicit description of stable spin curves. The formulas we obtain for the class of these divisors are identical over both $\ts_g^-$ and $\ts_g^+$, therefore we sometimes use the symbol $\ts_g^{\mp}$ (or even $\ts_g$), to denote one of the two spin moduli spaces.

We start with a  curve $[C]\in \cM_g$ and denote as usual by $Q_C:=M_{K_C}^{\vee}$ the associated \emph{Lazarsfeld bundle} \cite{L} defined via the exact sequence on $C$
$$0\rightarrow M_{K_C}\rightarrow H^0(C, K_C)\otimes \OO_C\stackrel{\mathrm{ev}}\rightarrow K_C\rightarrow 0.$$ Note that $Q_{C}$ is a semistable vector bundle on $C$ (even stable, when the curve $C$ is non-hyperelliptic), and $\mu(Q_{C})=2$. For integers $0\leq i\leq g-1$, one defines the \emph{divisorial difference variety}
$C_{g-i-1}-C_i\subset \mathrm{Pic}^{g-2i-1}(C)$ as being the image of
the difference map
$$\phi:C_{g-i-1}\times C_i \rightarrow \mbox{Pic}^{g-2i-1}(C), \mbox{ } \ \phi(D,
E):=\OO_C(D-E).$$ The main result from \cite{FMP} provides a scheme-theoretic identification of divisors on the Jacobian variety
\begin{equation}\label{mrc}
C_{g-i-1}-C_i=\Theta_{\wedge^i Q_C}\subset \mathrm{Pic}^{g-2i-1}(C),\end{equation}
where the right-hand-side denotes
the \emph{Raynaud locus} \cite{R} $$\Theta_{\wedge^i Q_C}:=\{\eta\in \mathrm{Pic}^{g-2i-1}(C): H^0(C, \wedge^i Q_C\otimes \eta)\neq 0\}.$$
The non-vanishing $H^0(C, \wedge^i Q_C\otimes \xi)\neq 0$ for all line bundles $\xi=\OO_C(D-E)$, where $D \in C_{g-i-1}$ and $E \in C_i$, follows from \cite{L}. The thrust of \cite{FMP}
is that the reverse inclusion $\Theta_{\wedge^i Q_C}\subset C_{g-i-1}-C_i$ also holds. Moreover, identification (\ref{mrc}) shows that, somewhat similarly to Riemann's Singularity Theorem,  the product $C_{g-i-1}\times C_i$ can be thought of as a canonical desingularization of the generalized theta-divisor $\Theta_{\wedge^i Q_C}$.
\vskip 5pt
We fix integers $r, s>0$ and set $d:=rs+r, g:=rs+s$, therefore the Brill-Noether number $\rho(g, r, d)=0$. We assume moreover that
$d\equiv 0 \mbox{ mod } 2$, that is, either $r$ is even or $s$ is odd, and write $d=2i$. We define the following locus in the spin moduli space
$\cS_g^{\mp}$:
$$\cU_{g, d}^r:=\{[C, \eta]\in \cS_g^{\mp}: \exists L\in W^r_d(C) \mbox{ such that } \eta\otimes L^{\vee}\in C_{g-i-1}-C_i\}.$$
Using (\ref{mrc}), the condition $[C, \eta]\in \cU_{g, d}^r$ can be rewritten in a determinantal way as, $$H^0(C, \wedge^i M_{K_C}\otimes \eta\otimes L)\neq 0.$$
Tensoring by $\eta\otimes L$ the exact sequence coming from the definition of $M_{K_C}$, namely
$$0\longrightarrow \wedge^i M_{K_C}\longrightarrow \wedge^i H^0(C, K_C)\otimes \OO_C\longrightarrow \wedge^{i-1} M_{K_C}\otimes K_C\longrightarrow 0,$$ then taking global sections and finally using that $M_{K_C}$ (hence all of its exterior powers)  are semi-stable vector bundles, we find that $[C, \eta]\in \cU_{g, d}^r$ if and only if the map
\begin{equation}\label{det}
\phi(C, \eta, L): \wedge^i H^0(C, K_C)\otimes H^0\bigl(C,  \eta\otimes L)\rightarrow H^0(C, \wedge^{i-1} M_{K_C}\otimes K_C\otimes \eta\otimes L\bigr)
\end{equation}
is not an isomorphism for a certain $L\in W^r_d(C)$. Since $\mu(\wedge^{i-1} M_{K_C}\otimes K_C\otimes \eta\otimes L)\geq 2g-1$ and $\wedge^{i-1} M_{K_C}$ is a semi-stable vector bundle on $C$, it follows that
$$h^0(C, \wedge^{i-1} M_{K_C}\otimes K_C\otimes \eta\otimes L)=\chi(C, \wedge^{i-1} M_{K_C}\otimes K_C\otimes \eta\otimes L)={g\choose i}d.$$
We assume that $h^1(C, \eta\otimes L)=0$. This condition is satisfied outside a locus of $\cS_g^{\mp}$ of codimension at least $2$; if $H^1(C,  \eta\otimes L)\neq 0$, then $H^1(C, K_C\otimes L^{\otimes (-2)})\neq 0$, in particular the Petri map $$\mu_0(C, L):H^0(C, L)\otimes H^0(C, K_C\otimes L^{\vee})\rightarrow H^0(C, K_C)$$ is not injective. Then $h^0(C, L\otimes \eta)=d$ and we note that $\phi(C, \eta, L)$ is a map between vector spaces of the same rank. This obviously suggests a determinantal presentation of $\cU_{g, d}^r$ as the (push-forward of) a degeneracy locus between vector bundles of the same rank.
In what follows we extend this presentation over a partial compactification of $\ts_g^{\mp}$. We refer to \cite{FL} Section 2 for a similar
calculation over the Prym moduli stack $\overline{\textbf{R}}_g$.
\vskip 4pt
We denote by $\textbf{M}_g^0\subset \textbf{M}_g$ the open substack
classifying curves $[C]\in \cM_g$ such that $W_{d-1}^r(C)=
\emptyset$, $W_d^{r+1}(C)= \emptyset$ and moreover $H^1(C, L\otimes \eta)=0$, for every $L\in W^r_d(C)$ and each odd-theta characteristic $\eta\in \mbox{Pic}^{g-1}(C)$. From general Brill-Noether theory one knows that
$\mbox{codim}(\cM_g-\cM_g^0, \cM_g)\geq 2$. Then we define $\widetilde{\Delta}_0\subset \Delta_0$ to be the open substack consisting of $1$-nodal
stable curves $[C_{yq}:=C/y\sim q]$, where $[C]\in \cM_{g-1}$ is a curve satisfying the Brill-Noether theorem and $y, q\in C$. We then set $\rem_g^0:=\textbf{M}_g^0\cup \widetilde{\Delta}_0$, hence $\rem_g^0\subset \pem_g$ and then $\ts_g^0:=\pi^{-1}(\rem_g^0)=(\ts_g^0)^+ \cup (\ts_g^0)^-$. Following
\cite{EH1}, \cite{F1}, we consider the  proper Deligne-Mumford stack
$$\sigma_0:\mathfrak G^r_d\rightarrow \rem_g^0$$
classifying pairs $[C, L]$ with $[C]\in \mm_g^0$ and $L\in W^r_d(C)$.
For any curve $[C]\in \mm_{g}^0$ and $L\in W^r_d(C)$, we have that
$h^0(C, L)=r+1$, that is, $\mathfrak{G}^r_d$ parameterizes only complete linear series.
For a point $[C_{yq}:=C/y\sim q]\in \widetilde{\Delta}_0$, we have the identification
$$\sigma_0^{-1}\bigl[C_{yq}\bigr]=\{L\in W^r_d(C): h^0(C,
L\otimes \OO_C(-y-q))=r\},$$
that is, we view linear series on singular curves as linear series on the normalization such that the divisor of the nodes imposes only one condition.  We denote by $f^r_d:\mathfrak C_{g,
d}^r:=\rem_{g, 1}^0\times_{\rem_g^0} \mathfrak G^r_d\rightarrow
\mathfrak G^r_d$ the pull-back of the universal curve $p:\rem_{g,
1}^0\rightarrow \rem_g^0$ to $\mathfrak G^r_{d}$. Once we have chosen a
Poincar\'e bundle $\L$ on $\mathfrak C^r_{g, d}$, we can form the
three codimension $1$ tautological classes in $A^1(\mathfrak
G^r_d)$:
\begin{equation}\label{tautological}
\mathfrak{a}:=(f^r_d)_*\bigl(c_1(\L)^2\bigr), \ \mathfrak{b}:=(f^r_d)_*\bigl(c_1(\L)\cdot
c_1(\omega_{f_d^r})\bigr), \mbox{ }
\mathfrak{c}:=(f^r_d)_*\bigl(c_1(\omega_{f_d^r})^2\bigr)=(\sigma_0)^*\bigl((\kappa_1)_{
\rem_g^0}\bigr).
\end{equation}
The dependence on $\mathfrak a, \mathfrak b, \mathfrak c$ on the choice of $\L$ is discussed in both \cite{F2} and \cite{FL}. We introduce
the stack of  $\mathfrak g^r_d$'s on spin curves
$$\sigma:\mathfrak G^r_d(\ts_g^0/\rem_g^0):=\ts_g^{0} \times_{\rem_g^0} \mathfrak
G^r_d\rightarrow \ts_g^0$$ and then the corresponding universal spin curve over the $\mathfrak g^r_d$ parameter space
$$f':\mathcal{C}_d^r:=\mathcal{C}\times _{\ts_g^0} \AUX\rightarrow \AUX.$$ We note that $f'$ is a family of quasi-stable curves carrying at the same time a spin structure as well as a $\mathfrak g^r_d$. Just like in \cite{FL},
the boundary divisors of $\AUX$ are denoted by the same symbols, that is, one sets
$A_0':=\sigma^*(A_0')$  and
$B_0':=\sigma^*(B_0')$ and then $$\alpha_0:=[A_0'], \ \ \beta_0:=[B_0']\in A^1(\AUX).$$
We observe that two tautological line bundles live on $\mathcal{C}^r_d$, namely the pull-back of the universal spin bundle
$\mathcal{P}_d^r\in \mathrm{Pic}(\cC_d^r)$  and a  Poincar\'e bundle
$\mathcal{L}\in \mbox{Pic}(\mathcal{C}_d^r)$ singling out the $\mathfrak g^r_d$'s, that is, $\mathcal{L}_{|f'^{-1}[X, \eta, \beta, L]}=L\in W^r_d(C),$
for each point $[X, \eta, \beta, L]\in \AUX$. Naturally, one also has the
classes $\mathfrak{a}, \mathfrak{b}, \mathfrak{c}\in A^1(\AUX)$ defined by the formulas
(\ref{tautological}).

The following result is easy to prove and we skip  details:
\begin{proposition}\label{tauto}
We denote by $f':\cC_{d}^r \rightarrow \AUX $ the universal quasi-stable spin curve and by $\mathcal{P}_d^r\in \mathrm{Pic}(\cC_d^r)$ the universal spin bundle of relative degree
$g-1$. One has the following formulas in $A^1(\AUX)$:
\begin{enumerate}
\item $f'_*(c_1(\omega_{f'})\cdot c_1(\mathcal{P}_d^r))=\frac{1}{2}\mathfrak c$.
\item $f'_*(c_1(\mathcal{P}_d^r)^2)=\frac{1}{4}\mathfrak c-\frac{1}{2} \beta_0.$
\item $f'_*(c_1(\mathcal{L})\cdot c_1(\mathcal{P}_d^r))=\frac{1}{2}\mathfrak b.$
\end{enumerate}
\end{proposition}

We determine the class of a compactification of $\cU_{g, d}^r$ by pushing-forward a codimension $1$ degeneracy locus via the map $\sigma:\AUX\rightarrow \ts_g^0$. To that end, we define a sequence of
tautological vector bundles on $\AUX$: First, for $l\geq 0$ we set $$\cA_{0, l}:=f'_*(\L\otimes
\omega_{f'}^{\otimes l} \otimes \P_d^r).$$ It is easy to verify that $R^1f'_*(\L\otimes
\omega_{f'}^{\otimes l} \otimes \P_d^r)=0$, hence $\cA_{0, l}$ is locally free
 over $\AUX$ of rank equal to $h^0(X, L\otimes \omega_X^{\otimes l}\otimes
 \eta)=l(2g-2)+d$.  Next we introduce the global
Lazarsfeld vector bundle $\cM$ over $\mathcal{C}_d^r$ by the exact
sequence
$$0\longrightarrow \cM\longrightarrow (f')^*\bigl(f'_*\omega_{f'}\bigr)\longrightarrow
\omega_{f'}\longrightarrow 0,$$ and then  for all integers $a, j\geq
1$ we define the sheaf over $\AUX$
$$\cA_{a, j}:=f'_*(\wedge^a \cM\otimes \omega_{f'}^{\otimes j}\otimes
\L \otimes \P^r_d).$$ In a way similar to \cite{FL} Proposition 2.5 one shows
 that $R^1 f'_*\bigl(\wedge^a \cM\otimes \omega_{f'}^{\otimes (i-a)}\otimes
\L\otimes \P^r_d\bigr)=0$,
therefore by Grauert's theorem $\cA_{a, i-a}$ is a vector bundle over $\AUX$ of rank
$$\mathrm{rk}(\cA_{a, i-a})=\chi\bigl(X, \wedge^a M_{\omega_X}\otimes \omega_X^{\otimes
(i-a)}\otimes L \otimes \eta\bigr)=2(i-a)g{g-1\choose a}.$$
Furthermore, for all  $1\leq a\leq i-1$,
the vector bundles $\cA_{a, i-a}$ sit in exact sequences
\begin{equation}\label{recursion}
0\longrightarrow \cA_{a, i-a}\longrightarrow \wedge^a f'_*(\omega_{f'})\otimes
\cA_{0, i-a}\longrightarrow \cA_{a-1, i-a+1}\longrightarrow 0,
\end{equation}
where the right exactness boils down to showing that $H^1(X, \wedge^a M_{\omega_X}\otimes \omega_X^{\otimes (i-a)}\otimes
\eta\otimes L)=0$ for all $[X, \eta, \beta, L]\in \AUX$.

We denote as usual $\mathbb E:=f'_*(\omega_{f'})$ the Hodge bundle over $\AUX$ and then note that there exists a vector bundle map
\begin{equation}\label{phi}
\phi: \wedge^i \mathbb E\otimes \cA_{0, 0}\rightarrow \cA_{i-1, 1}
\end{equation}
between vector bundles of the same rank over $\AUX$. For $[C, \eta, L]\in \sigma^{-1}(\cM_g^0)$ the fibre of this morphism is precisely the map $\phi(C, \eta, L)$ defined by (\ref{det}).

\begin{theorem} The vector bundle morphism $\phi:\wedge^i \mathbb E\otimes \cA_{0, 0}\rightarrow \cA_{i-1, 1}$ is generically non-degenerate
over $\AUX$. It follows that $\cU_{g, d}^r$ is an effective divisor over $\cS_g^+$ for all $s\geq 1$, and over $\cS_g^-$ as well for
$s\geq 2$.
\end{theorem}
\begin{proof} We specialize $C$ to a hyperelliptic curve, and denote by $A\in W^1_2(C)$ the hyperelliptic involution.
The Lazarsfeld bundle splits into a sum of line bundles $Q_C\cong A^{\oplus (g-1)}$, therefore the condition $H^0(C, \wedge^i M_{K_C}\otimes\eta\otimes L)=0$ translates into
$H^0(C, \eta\otimes A^{\otimes i}\otimes L^{\vee})=0$. Suppose that $h^0(C, \eta\otimes A^{\otimes i}\otimes L^{\vee})\geq 1$ for any $L=A^{\otimes r}\otimes \OO_C(x_1+\cdots+x_{d-2r})\in W^r_d(C)$, where
the $x_1, \ldots, x_{d-2r}\in C$ are arbitrarily chosen points. This implies that $h^0(C, \eta\otimes A^{\otimes (i-r)})\geq d-2r+1$.  Any theta-characteristic on $C$ is of
the form $$\eta=A^{\otimes m}\otimes \OO_C(p_1+\cdots +p_{g-2m-1}),$$ where $1\leq m\leq (g-1)/2$ and $p_1, \ldots, p_{g-2m-1}\in C$ are Weierstrass points. Choosing a theta-characteristic on $C$ for which $m\leq i-r-1$ (which can be done in all cases except on $\cS_g^-$ when $i=r$), we obtain that
$h^0(C, \eta\otimes A^{\otimes (i-r)})\leq d-2r$, a contradiction.
\end{proof}

\noindent
\emph{Proof of Theorem \ref{spinclass}.} To compute the class of the degeneracy locus of $\phi$ we use
repeatedly the exact sequence (\ref{recursion}). We
write the following identities in $A^1(\AUX)$:
$$c_1\bigl(\cA_{i-1, 1}-\wedge^i \mathbb E \otimes \cA_{0,
0}\bigr)=\sum_{l=0}^i (-1)^{l-1} c_1(\wedge^{i-l} \mathbb E \otimes \cA_{0,
l})=$$
$$=\sum_{l=0}^i (-1)^{l+1} \Bigl((2l(g-1)+d){g-1\choose
i-l-1}c_1(\mathbb E)+{g\choose i-l}c_1(\cA_{0, l})\Bigr).$$
Using Proposition \ref{tauto} one can show via the Grothendieck-Riemann-Roch formula applied to $f':\cC_{d}^r\rightarrow \AUX$ that one has that
$$c_1(\cA_{0, l})=\lambda+\Bigl(\frac{l^2}{2}-\frac{1}{8}\Bigr)\mathfrak{c}+\frac{1}{2}\mathfrak a+l\mathfrak{b}-\frac{1}{4}\beta_0\in A^1(\AUX).$$
To determine $\sigma_*\bigl(c_1(\cA_{i-1, 1}-\wedge^i \mathbb E)\bigr)\in A^1(\ts_g)$ we use \cite{F1}, \cite{Kh}: If $$N:=\mbox{deg}(\sigma)=\#(W^r_d(C))$$ denotes the number of $\mathfrak g^r_d$'s on a general curve $[C]\in \cM_g$, then there exists a precisely described choice of a Poincar\'e bundle on $\mathfrak C_{g, d}^r$ such that the push-forwards
of the tautological classes on $\AUX$ are given as follows (cf. \cite{F1}, \cite{Kh} and especially \cite{FL} Section 2, for the similar argument in the Prym case):
$$\sigma_*(\mathfrak a)=\frac{dN}{(g-1)(g-2)}\Bigl((gd-2g^2+8d-8g+4)\lambda+\frac{1}{6}(2g^2-gd+3g-4d-2)(\alpha_0+2\beta_0)\Bigr)$$
and
$$\sigma_*(\mathfrak b)=\frac{dN}{2g-2}\Bigl(12\lambda-\alpha_0-2\beta_0\Bigr)\in A^1(\AUX).$$
One notes that $c_1(\cA_{i-1, 1}-\wedge^i \mathbb E\otimes \cA_{0, 0})\in A^1(\AUX)$ does not depend of the Poincar\'e bundle. Using the previous formulas, after some arithmetic, one computes the class of the partial compactification of $\cU_{g, d}^r$ and finishes the proof.
\hfill $\Box$

When $s=2a+1$, hence $g=(2a+1)(r+1)$ and $d=2r(a+1)$, our calculation shows that $$\overline{\cU}_{g, d}^r\equiv c_{a, r}\bigl(\bar{\lambda}\ \lambda-\bar{\alpha_0}\ \alpha_0-\bar{\beta}_0\ \beta_0-\cdots)\in \mbox{Pic}(\ss_g^{\mp}),$$ where $c_{a, r}\in \mathbb Q_{>0}$ is explicitly known and
{\small{
$$\bar{\lambda}=12r^3-12r^2-48a^2+96a^3+48r^4a+2208r^3a^3+1968r^3a^2+3936r^2a^3+2208ra^3+
552r^3a+3984r^2a^2+
$$$$1080r^2a+2160ra^2+528ra+192r^4a^4+384r^4a^3+768r^3a^4+
960r^2a^4+240r^4a^2+384ra^4,$$

$$\bar{\alpha}_0=220ra^2+536r^2a^3+32r^4a^4+36ra+24a^3+328r^3a^3+296ra^3+8r^4a+64r^4a^3+3r^3
+468r^2a^2+$$
$$128r^3a^4+74r^3a+40r^4a^2+160r^2a^4+64ra^4+268r^3a^2+110r^2a-3r^2-12a^2$$
and
$$\bar{\beta}_0=96ra+64r^4a^4+16r^4a+416ra^2+928r^2a^3+448ra^3+208r^2a+608r^3a^3+256r^3a^4+112r^3a+$$
$$80r^4a^2+320r^2a^4+128ra^4+464r^3a^2+128r^4a^3+816r^2a^2.$$
}
}
\vskip 4pt
These formulas, though unwieldy, carry a great deal of information about $\ss_g$. In the simplest case,  $s=1$ (that is, $a=0$) and $r=g-1$, then necessarily $L=K_C\in W_{2g-2}^{g-1}(C)$ and the condition $\eta-K_C\in -C_{g-1}$ is equivalent to $H^0(C, \eta)\neq 0$. In this way we recover  the theta-null divisor $\thet$ on $\ss_g^+$, or more precisely also taking into account multiplicities \cite{F2},  $$\cU_{g, 2g-2}^{g-1}=2\cdot \Theta_{\mathrm{null}}.$$ At the same time,  on
$\cS_g^+$ one does not get a divisor at all. In particular, we find that
$$\overline{\mathcal{U}}_{g, 2g-2}^{g-1}\equiv 2\cdot \overline{\Theta}_{\mathrm{null}}\equiv \frac{1}{2}\lambda-\frac{1}{8}\alpha_0-0\cdot \beta_0-\cdots \in \mathrm{Pic}(\ss_g^+).$$
Another interesting case is when $r=2$, hence $g=3s, L\in W^2_{2s+2}(C)$ and the condition $\eta\otimes L^{\vee}\in C_{2s-2}-C_{s+1}$ is equivalent to requiring
that the embedded curve $C\stackrel{|\eta\otimes L|}\longrightarrow \PP^{2s+1}$ has an $(s+1)$-secant $(s-1)$-plane:
\begin{theorem}
For $g=3s, d=2s+2$, the class of the closure in $\overline{\textbf{S}}_g^{\mp}$ of the effective divisor
$$\cU_{g, d}^2:=\{[C, \eta]\in \cS_{3s}^{\mp}:\exists L\in W^2_{2s+2}(C) \ \mbox{ such that } \eta\otimes L^{\vee}\in C_{2s-2}-C_{s+1}\}$$
is given by the formula in $\mathrm{Pic}(\ss_g^{\mp})$:
$$\overline{\cU}_{g, d}^2\equiv{g\choose s+2}{g\choose s, s, s}\frac{1}{24g(g-1)^2(g-2)(s+1)^2}
\Bigl(4(216s^4+513s^3-348s^2-387s+18\bigr)\lambda-$$
$$-\bigl(144s^4+225s^3-268s^2-99s+10\bigr)\alpha_0
-\bigl(288s^4+288s^3+320s^2+32\bigr)\beta_0-\cdots\Bigr).
$$
\end{theorem}
For instance, for $g=9$, we obtain the class of the closure of the locus spin curves $[C, \eta]\in \cS_9^{\mp}$, for which
there exists a net $L\in W^2_8(C)$ such that $\eta\otimes L^{\vee}\in C_4-C_4$:
$$\overline{\cU}_{9, 8}^2\equiv 235\cdot 35 \Bigl(\frac{36}{5}\lambda-\alpha_0-\frac{428}{235}\beta_0-\cdots \Bigr)\in \mathrm{Pic}(\ss_9^{\mp}).$$

\section{The class of $\overline{\Theta}_{\mathrm{null}}$ on $\ss_g^+$: An alternative
proof using the Hurwitz stack}

We present an alternative way of computing the class of the divisor
$[\overline{\Theta}_{\mathrm{null}}]$ \ (in even genus), as the push-forward of a
determinantal cycle on a Hurwitz scheme of degree $k$ coverings of genus $g$ curves. We set
$$g=2k-2, \  r=1, \ d=k, $$
hence $\rho(g, 1, k)=0$, and use the notation from the previous section.
In particular, we have the proper morphism  $\sigma_0:\mathfrak
G^1_k\rightarrow \rem_g^0$ from the Hurwitz stack of $\mathfrak g^1_k$'s, and the universal spin curve over the Hurwitz stack

$$f':\mathcal{C}_1^k:=\mathcal{C}\times _{\ts_g^0} \aux
\rightarrow  \aux.$$
Once more, we introduce a number of vector bundles over $\aux$: First, we set
$\H:=f'_*(\L).$ By Grauert's theorem,  $\H$ is a vector bundle of
rank $2$ over $\aux$, having fibre $\H[X, \eta, \beta, L]=H^0(X, L)$, where $L\in W^1_k(X)$. Then for
$j\geq 1$ we define
$$\cB_{j}:=f'_*(\L^{\otimes j}\otimes
\P_k^1).$$ Since $R^1f'_*(\L^{\otimes j}\otimes \P_k^1)=0$, we find
that $\cB_{j}$ is a vector bundle over $\aux$ of rank equal to
$h^0(X, L^{\otimes j}\otimes \eta)=kj$.

\begin{proposition}\label{rr}
If $\mathfrak a, \mathfrak b, \mathfrak c$ are the codimension $1$
tautological classes on $\aux$ defined by (\ref{tautological}), then for all  $j\geq 1$ one has the
following formula in $A^1(\aux)$:
$$c_1(\cB_{j})=\lambda-\frac{1}{8}\mathfrak c +\frac{j^2}{2}\mathfrak a-\frac{j}{2}\mathfrak{b}
-\frac{1}{4}\beta_0.$$
\end{proposition}
\begin{proof}
We apply Grothendieck-Riemann-Roch to the morphism
$f':\mathcal{C}_k^1\rightarrow \aux$:
$$c_1(\cB_{j})=c_1\bigl(f'_{!}(\L^{\otimes j}\otimes \P^1_k)\bigr)=
$$
$$=f'_*\Bigl[\Bigl(1+c_1(\L^{\otimes j}\otimes
\P^1_k)+\frac{c_1^2(\L^{\otimes j} \otimes
\P^1_k)}{2}\Bigr)\Bigl(1-\frac{c_1(\omega_{f'})}{2}+\frac{c_1^2(\omega_{f'})
+[\mathrm{Sing}(f')]}{12}\Bigr)\Bigr]_2,$$ where
$\mathrm{Sing}(f')\subset \mathcal{X}_{k}^1$ denotes the codimension
$2$ singular locus of the morphism $f'$, therefore
$f'_*[\mathrm{Sing}(f')]=\alpha_0+2\beta_0$. We then use Mumford's
formula \cite{HM} pulled back from $\rem_g^0$ to $\aux$, to write that
$$\kappa_1=f'_*(c_1^2(\omega_{f'}))=12\lambda-(\alpha_0+2\beta_0)$$
and then note that $f'_*(c_1(\L)\cdot c_1(\mathcal{P}_k^1))=0$ (the
restriction of $\L$ to the exceptional divisor of
$f':\mathcal{C}_k^1\rightarrow \aux$ is trivial). Similarly, we note
that $f'_*(c_1(\omega_{f'})\cdot c_1(\mathcal{P}_k^1))=\mathfrak
c/2$. Finally, we write that
$f'_*(c_1^2(\mathcal{P}_k^1))=\mathfrak{c}/4-\beta_0/2$.
\end{proof}

For $j\geq 1$ there are natural vector bundle morphisms over $\aux$
$$\chi_j:\H\otimes \cB_j\rightarrow \cB_{j+1}.$$
Over a point $[C,
\eta_C, L]\in \cS_g^+\times_{\cM_g} \mathfrak G^1_k$ corresponding to an even theta-characteristic
$\eta_C$ and a pencil $L\in W_k^1(C)$, the morphism $\chi_j$
is given by multiplications of global sections
$$\chi_j[C, \eta, L]: H^0(C, L)\otimes H^0(C, L^{\otimes j}\otimes
\eta_C)\rightarrow   H^0(C, L^{\otimes (j+1)}\otimes \eta_C).$$ In
particular, $\chi_1:\H\otimes \cB_1\rightarrow \cB_2$ is a morphism
between vector bundles of the same rank.
 From the base point free pencil trick, the degeneration locus $Z_1(\chi_1)$ is (set-theoretically) equal to the inverse image
 $\sigma^{-1}\bigl(\overline{\Theta}_{\mathrm{null}}\cap (\ss_g^0)^+\bigr)$.

\begin{theorem}
We fix $g=2k-2$. The vector bundle morphism $\chi_1:\H\otimes
\cB_1\rightarrow \cB_2$ defined over $\aux$  is generically
non-degenerate and we have the following formula for the class of its degeneracy
locus:
$$[Z_1(\chi_1)]=c_1(\cB_2-\H\otimes \cB_1)=\frac{1}{2}\lambda-\frac{1}{8}\alpha_0+\mathfrak a-kc_1(\H)\in A^1(\aux).$$
The class of the push-forward $\sigma_*[Z_1(\chi_1)]$ to $\ss_g^+$ is given by the formula:
$$\sigma_*\bigl(c_1(\cB_2-\H\otimes \cB_1)\bigr)\equiv \frac{(2k-2)!}{k! (k-1)!}\Bigl(\frac{1}{2}\lambda-\frac{1}{8}\alpha_0-0\cdot \beta_0\Bigr)
\equiv \frac{2(2k-2)!}{k! (k-1)!}\ \overline{\Theta}_{\mathrm{null } \ |
\ss_g^+} \in \mathrm{Pic}(\ss_g^+).$$

\end{theorem}
\begin{proof} The first part follows directly from Theorem \ref{rr}.
 To determine the push-forward of codimension $1$ tautological classes to $(\ss_g^0)^+$, we use again \cite{F1}, \cite{Kh}:
One writes the following relations in $A^1((\ts_g^0)^+)=A^1((\ss_g^0)^+)$:
$$\sigma_*(\mathfrak a)=\mathrm{deg}(\mathfrak G^1_k/\rem_g^0)\ \Bigl(-\frac{3k(k+1)}{2k-3}\ \lambda+\frac{k^2}{2(2k-3)}(\alpha_0+2\beta_0)\Bigr),$$
$$ \sigma_*(\mathfrak b)=\mathrm{deg}(\mathfrak G^1_k/\rem_g^0)\ \Bigl(\frac{6k}{2k-3}\ \lambda-\frac{k}{2(2k-3)}(\alpha_0+2\beta_0)\Bigr),$$
and $$\sigma_*\bigl(c_1(\H))=\mathrm{deg}(\mathfrak
G^1_k/\rem_g^0)\ \Bigl(-3\frac{k+1}{2k-3}\
\lambda+\frac{k}{2(2k-3)}(\alpha_0+2\beta_0)\Bigr),$$ where $$N:=
\mathrm{deg}(\mathfrak G^1_k/\rem_g^0)=\frac{(2k-2)!}{k! (k-1)!}$$ denotes the
\emph{Catalan number} of linear series $\mathfrak g^1_k$ on a general curve of genus $2k-2$. This yields yet another proof of the main result from
\cite{F2}, in the sense that we compute the class of the divisor $\thet$ of vanishing theta-nulls:
$$\sigma_*(Z_1(\chi_1))=\mathrm{deg}(\mathfrak G^1_k/\rem_g^0)\bigl(\frac{1}{2}\lambda-\frac{1}{8}\alpha_0 \bigr)\equiv 2\mathrm{deg}(\mathfrak G^1_k/\rem_g^0)\ [\overline{\Theta}_{\mathrm{null} \ | (\ss_g^0)^+}].$$
\end{proof}
\begin{remark} The multiplicity $2$ appearing in the expression of
$\sigma_*(Z_1(\chi_1))$ is justified by the fact that $\mbox{dim }
\mbox{Ker}(\chi_1(t))=h^0(C, \eta)$ for every $[C, \eta, L]\in
\sigma^{-1}((\cS_g^0)^+)$.  This of course is always an even number. Thus we have the equality cycles
$$Z_1(\chi_1)=Z_2(\chi_1)=\{t\in \aux: \mbox{co-rank}( \phi_1(t))\geq
2\},$$ that is $\chi_1$ degenerates in codimension $1$ with corank
$2$, and $Z_1(\chi_1)$ is an everywhere non-reduced scheme.
\end{remark}

\section{The divisor of points of odd theta-characteristics}
In this section we compute the class of the divisor $\tet_{g, 1}$. The study of geometric divisors on $\mm_{g, 1}$ begins with \cite{Cu}, where the locus of Weierstrass points is determined:
$$\overline{\W_g} \equiv -\lambda+{g+1\choose 2}\psi-\sum_{i=1}^{g-1}{g-i+1\choose 2} \delta_{i: 1}\in \mathrm{Pic}(\mm_{g, 1}).$$
More generally, if $\bar{\alpha}:0\leq \alpha_0\leq \ldots \leq \alpha_r\leq d-r$ is a \emph{Schubert index} of type $(r, d)$ such that $\rho(g, r, d)-\sum_{i=0}^{r}\alpha_i=-1$,  one defines the \emph{pointed Brill-Noether divisor} $\cM_{g, d}^r(\bar{\alpha})$ as being the locus of pointed curves $[C, q]\in \cM_{g, 1}$ possessing
a linear series $l\in G^r_d(C)$ with ramification sequence $\alpha^l(q)\geq \bar{\alpha}$. It follows from \cite{EH3} that the cone spanned by the pointed Brill-Noether divisors on $\mm_{g, 1}$ is $2$-dimensional, with generators $[\overline{\W}_g]$ and the pull-back of the Brill-Noether class from
$\mm_g$. Our aim is to analyze the divisor $\tet_{g, 1}$, whose definition is arguably simpler than that of the divisors $\mm_{g, d}^r(\bar{\alpha})$, and which seems to have been overlooked until now. A consequence of the calculation is that (as expected) $[\tet_{g, 1}]$ lies outside the Brill-Noether cone of $\mm_{g, 1}$.

We begin by recalling basic facts about divisors on $\mm_{g, 1}$.  For $i=1, \ldots, g-1$, the divisor $\Delta_i$ on $\mm_{g, 1}$ is the closure of the locus of pointed curves $[C\cup D, q]$, where $C$ and $D$ are smooth curves of genus $i$ and $g-i$ respectively, and $q\in C$.
Similarly, $\Delta_{\mathrm{irr}}$ denotes the closure in $\mm_{g, 1}$ of the locus of irreducible $1$-pointed stable curves. We set $\delta_i:=[\Delta_i], \delta_{\mathrm{irr}}:=[\Delta_{\mathrm{irr}}]\in \mbox{Pic}(\mm_{g, 1})$, and recall that $\psi \in \mbox{Pic}(\mm_{g, 1})$ is the universal cotangent class. Clearly, $p^*(\delta_{\mathrm{irr}})=\delta_{\mathrm{irr}}$ and $p^*(\delta_i)=\delta_i+\delta_{g-i}\in \mathrm{Pic}(\mm_{g, 1})$ for $1\leq i\leq [g/2]$. For $g\geq 3$, the group $\mbox{Pic}(\mm_{g, 1})$ is freely generated by the classes $\lambda, \psi,
\delta_{\mathrm{irr}}, \delta_1, \ldots, \delta_{g-1}$, cf. \cite{AC1}. When $g=2$, the same classes generate $\mbox{Pic}(\mm_{2, 1})$ subject to the \emph{Mumford relation} $$\lambda=\frac{1}{10}\delta_{\mathrm{irr}}+\frac{1}{5}\delta_1,$$ expressing that $\lambda$ is a boundary class. We expand the class $[\tet_{g, 1}]$ in this basis of $\mbox{Pic}(\mm_{g, 1})$,
$$\tet_{g, 1} \equiv a\lambda+b\psi-b_{\mathrm{irr}}\delta_{\mathrm{irr}}-\sum_{i=1}^{g-1} b_i\delta_i\in \mathrm{Pic}(\mm_{g, 1}),$$
 and determine the coefficients in a classical way, by understanding the restriction of  $\tet_{g, 1}$ to sufficiently many geometric subvarieties of $\mm_{g, 1}$.
To ease calculations, we set
$$N_g^-:=2^{g-1}(2^g-1) \ \mbox{  and } \ N_g^+:=2^{g-1}(2^g+1),$$
to be the number of odd (respectively even) theta-characteristic on a curve of genus $g$.

We define some test-curves in the boundary of $\mm_{g, 1}$. For an integer $2\leq i\leq g-1$, we choose general (pointed) curves $[C]\in \cM_{i}$ and $[D, x, q]\in \cM_{g-i, 2}$. In particular, we may assume that $x, q\in D$ do not appear in the support of any odd theta-characteristic $\eta_D^-$ on $D$, and that $h^0(D, \eta_D^+)= 0$, for any even theta-characteristic $\eta_D^+$. By joining $C$ and $D$ at a variable point $x\in C$, we obtain a family of $1$-pointed stable curves
$$F_{g-i}:=\{[C\cup_x D, \  q]: x\in C\} \subset \Delta_{g-i}\subset \mm_{g, 1},$$
where the marked point $q\in D$ is fixed.
It is clear that $F_{g-i}\cdot \delta_{g-i}=2-2i$, \ $F_{g-i}\cdot \lambda=F_{g-i}\cdot \psi=0$. Moreover, $F_{g-i}$ is disjoint from all the other boundary divisors of $\mm_{g, 1}$.

\begin{proposition}\label{fi}
For each $2\leq i\leq g-1$, one has that $b_{g-i}=N_i^-\cdot N_{g-i}^+/2$.
\end{proposition}
\begin{proof} We observe that the curve $F_{g-i}\times_{\mm_{g, 1}} \ss_g^{-}$ splits into $N_{i}^+\cdot N_{g-i}^- + N_i^-\cdot N_{g-i}^{+}$  irreducible components, each isomorphic to $C$,  corresponding to a choice of a pair of theta-characteristics of opposite parities on $C$ and $D$ respectively.
Let $t\in F_{g-i}\cdot \tet_{g, 1}$ be an arbitrary point in the intersection, with underlying stable curve $C\cup_x D$, and spin curves $\bigl([C, \eta_C], [D, \eta_D]\bigr)\in \cS_{i}\times \cS_{g-i}^{}$ on the two components.

Suppose first that $\eta_C=\eta_C^+$ and $\eta_D=\eta_D^-$, that is, $t$ corresponds to an even theta-characteristic on $C$ and an odd theta-characteristic on $D$. Then there exist non-zero sections $\sigma_C\in H^0\bigl(C, \eta_C^+\otimes \OO_C((g-i)x)\bigr)$ and $\sigma_D\in H^0\bigl(D, \eta_D^-\otimes \OO_D(ix)\bigr)$ such that
\begin{equation}\label{comp}
 \mathrm{ord}_x(\sigma_C)+\mathrm{ord}_x(\sigma_D)\geq g-1, \ \mbox{ and } \sigma_D(q)=0.
\end{equation}
In other words, $\sigma_C$ and $\sigma_D$ are the aspects of a limit $\mathfrak g^0_{g-1}$ on $C\cup_x D$ which vanishes at $q\in D$.
Clearly, $\mbox{ord}_x(\sigma_C)\leq g-i-1$, hence $\mbox{div}(\sigma_D)\geq ix+q$, that is, $q\in \mbox{supp}(\eta_D^-)$. This contradicts the generality assumption on $q\in D$, so this situation does not occur.
\vskip 3pt

Thus, we are left to consider the case $\eta_C=\eta_C^-$ and $\eta_D=\eta_D^+$. We denote again by $\sigma_C\in H^0(C, \eta_C^-\otimes \OO_C((g-i)x))$ and $\sigma_D\in H^0(D, \eta_D^+\otimes \OO_D(ix))$ the sections satisfying the compatibility relations (\ref{comp}). The condition $h^0(D, \eta_D^+\otimes \OO_D(x-q))\geq 1$ defines a correspondence on $D\times D$, cf. \cite{DK}, in particular, we can choose the points $x, q\in D$ general enough such that $H^0(D, \eta_D^+\otimes \OO_D(x-q))=0$. Then $\mbox{ord}_x(\sigma_D)\leq i-2$, thus $\mbox{ord}_x(\sigma_C)\geq g-i+1$. It follows that we must have equality
$\mbox{ord}_x(\sigma_C)=g-i+1$, and then, $x\in \mbox{supp}(\eta_C^-)$. An argument along the lines of \cite{EH3} Lemma 3.4, shows that each of these intersection points has to be counted with multiplicity $1$, thus $F_{g-i}\cdot \tet_{g, 1}=\#\mathrm{supp}(\eta_C^-)\cdot N_{i}^{-} \cdot N_{g-i}^+$.
We conclude by noting that $(2i-2)b_{g-i}=F_{g-i}\cdot \tet_{g, 1}$.
\end{proof}

\begin{proposition}\label{psicoeff}
The relation $b=N_g^{-} /2$ \ holds.
\end{proposition}
\begin{proof} Having fixed a general curve $[C]\in \mm_{g}$, by considering the fibre $p^*([C])$ inside the universal curve, one writes the identity $(2g-2)b=p^*([C])\cdot \tet_{g, 1}=(g-1) N_g ^-$.
\end{proof}
We compute the class of the restriction of  the divisor $\Theta_{g, 1}$ over $\cM_{g, 1}$:
\begin{proposition}\label{lambdacoeff}
One has the equivalence \  $\Theta_{g, 1} \equiv N_g^{-} (\psi/2+\lambda/4)\in \mathrm{Pic}(\cM_{g, 1})$.
\end{proposition}
\begin{proof} We consider the universal pointed spin curve \ $\mbox{pr}:\textbf{S}_{g, 1}^-:=\textbf{S}_g^{-}\times_{\textbf{M}_g} \textbf{M}_{g, 1}\rightarrow \textbf{M}_{g, 1}$. As usual, $\P\in \mbox{Pic}(\textbf{S}_{g, 1}^-)$ denotes the universal spin bundle, which over the stack $\textbf{S}_{g, 1}^-$, is a root of the dualizing sheaf $\omega_{\mathrm{pr}}$, that is, $2c_1(\P)=\mathrm{pr}^*(\psi)$. We introduce the divisor
$$\cZ:=\{[C, \eta, q]\in \cS_{g, 1}^-: q\in \mathrm{supp}(\eta)\}\subset \cS_{g, 1}^-,$$ and clearly $\Theta_{g, 1}:=\mathrm{pr}_*(\cZ)$. We write $[\cZ]=c_1(\P)-c_1\bigl(\mathrm{pr}^*(\mathrm{pr}_*(\P))\bigr)$, and take into account that $c_1(\mathrm{pr}_!(\P))=2c_1(\mathrm{pr}_*(\P))=-\lambda/2$. The rest follows by applying the projection formula.
\end{proof}
In order to determine the remaining coefficients $b_0, b_1$, we study the pull-back of $\tet_{g, 1}$ under the map $\nu: \mm_{1, 2}\rightarrow \mm_{g, 1}$,
given by $\nu([E, x, q]):=[C\cup_x E, q]\in \mm_{g, 1}$,
where $[C, x]\in \cM_{g-1, 1}$ is a fixed general pointed curve.

On the surface $\mm_{1, 2}$, if we denote a general element by $[E, x, q]$, one has the following relations between divisors classes, see \cite{AC2}:
$$\psi_x=\psi_q,\ \lambda=\psi_x-\delta_{0: xq},\ \delta_{\mathrm{irr}}=12(\psi_x-\delta_{0: xq}).$$  We describe the pull-back map  $\nu^*: \mbox{Pic}(\mm_{g, 1})\rightarrow \mbox{Pic}(\mm_{1, 2})$ at the level of divisors:
$$\nu^*(\lambda)=\lambda,\ \ \nu^*(\psi)=\psi_q, \ \ \nu^*(\delta_{\mathrm{irr}})=\delta_{\mathrm{irr}}, \ \ \nu^*(\delta_1)=-\psi_x, \ \ \nu^*(\delta_{g-1})=\delta_{0: xq}.$$
By direct calculation, we write $\nu^*(\tet_{g, 1})\equiv (a+b-12b_0+b_1)\psi_x-(a+b_{g-1}-12b_0)\delta_{0: xq}$.
We compute $b_0$ and $b_1$ by describing $\nu^*(\tet_{g, 1})$ viewed as an explicit divisor on $\mm_{1, 2}$:
\begin{proposition}\label{nustar}
One has the relation $\nu^*(\tet_{g, 1})\equiv N_{g-1}^- \cdot \overline{\mathfrak T}_2\in \mathrm{Pic}(\mm_{2, 1})$, \ where
$$\mathfrak{T}_2:=\{[E, x, q]\in \cM_{1, 2}: 2x\equiv 2q\}.$$
\end{proposition}
\begin{proof} We fix an arbitrary point $t:=[C\cup_x E, q]\in \nu^*(\tet_{g, 1})$. Suppose first that $E$ is a smooth elliptic curve, that is, $j(E)\neq \infty$ and $x\neq q$. Then there exist theta-characteristics of opposite parities $\eta_C, \eta_E$ on $C$ and $E$ respectively, together with non-zero sections $$\sigma_C\in H^0\bigl(C, \eta_C\otimes \OO_C(x)\bigr)\ \mbox{ and } \  \sigma_E\in H^0\bigl(E, \eta_E\otimes \OO_E((g-1)x)\bigr),$$ such that $\sigma_E(q)=0$ and \
$\mbox{ord}_x(\sigma_C)+\mbox{ord}_x(\sigma_E)\geq g-1.$

First we assume that $\eta_C=\eta_C^+$ and $\eta_E=\eta_E^-$, thus, $\eta_E=\OO_E$. Since $H^0(C, \eta_C^+)=0$, one obtains that
$\mbox{ord}_x(\sigma_C)=0$, that is $\mbox{ord}_{x}(\sigma_E)=g-1$, which is impossible, because $\sigma_E$ must vanish at $q$ as well. Thus, one is lead to study the remaining case, when
$\eta_C=\eta_C^-$ and $\eta_E=\eta_E^+$. Since $x\notin \mbox{supp}(\eta_C^-)$, we obtain $\mbox{ord}_x(\sigma_C)\leq 1$, and then by compatibility, the last inequality becomes equality, while $\mbox{ord}_x(\sigma_E)=g-2$, hence $\eta_E^+=\OO_E(x-q)$, or equivalently, $[E, x, q]\in \mathfrak{T}_2$.
The multiplicity $N_{g-1}^-$ in the expression of $\nu^*(\tet_{g, 1})$ comes from the choices for the theta-characteristics $\eta_C^-$, responsible for the $C$-aspect of a limit $\mathfrak{g}_{g-1}^0$ on $C\cup_x E$. It is an easy moduli count to show that the cases when $j(E)=\infty$, or  $[E, x, q]\in \delta_{0: xq}$ (corresponding to the situation when $x$ and $q$ coalesce on $E$), do not occur generically on a component of $\nu^*(\tet_{g, 1})$.
\end{proof}
\begin{proposition}\label{class2tors}
$\overline{\mathfrak{T}}_2$ is an irreducible divisor on $\mm_{1, 2}$ of class \ $\overline{\mathfrak{T}}_{2}\equiv 3\psi_x\in \mathrm{Pic}(\mm_{1, 2})$.
\end{proposition}
\begin{proof} We write $\overline{\mathfrak{T}}_2\equiv \alpha\psi_x-\beta \delta_{0: xq}\in \mathrm{Pic}(\mm_{1, 2})$, and we need to understand the intersection of $\overline{\mathfrak{T}}_2$ with two test curves in $\mm_{1, 2}$. First, we fix a general point $[E, q]\in \mm_{1, 1}$ and consider the family
$E_1:=\{[E, x, q]: x\in E\}\subset \mm_{1, 2}$. Clearly, $E_1\cdot \delta_{0: xq}=E_1\cdot \psi_x=1$. On the other hand $E_1\cdot \overline{\mathfrak{T}}_2$ is a $0$-cycle simply supported at the points $x\in E-\{q\}$ such that $x-q\in \mbox{Pic}^0(E)[2]$, that is, $E_1\cdot \overline{\mathfrak{T}}_2=3$. This yields the relation $\alpha-\beta=3$.
\vskip 6pt
As a second test curve, we denote by $[L, u, x, q]\in \mm_{0, 3}$ the rational $3$-pointed rational curve, and define the pencil $R:=\{[L\cup_u E_{\lambda}, x, q]: \lambda\in \PP^1\}\subset \mm_{1, 2}$, where $\{E_{\lambda}\}_{\lambda \in \PP^1}$ is a pencil of plane cubic curves. Then $R\cap \overline{\mathfrak{T}}_2=\emptyset$. Since $R\cdot \lambda=1$ and $R\cdot \delta_{\mathrm{irr}}=12$, we obtain the additional relation $\beta=0$, which completes the proof.
\end{proof}

Putting together Propositions \ref{fi}, \ref{lambdacoeff} and \ref{class2tors}, we obtain the system of equations
$$a+b_{g-1}-12b_{\mathrm{irr}}=0, \ a-12b_{\mathrm{irr}}+b+b_1=3N_{g-1}^-,\ a=\frac{1}{4} N_g^-, \ b=\frac{1}{2} N_g^-, \ b_1=\frac{3}{2}N_{g-1}^-.$$
Thus $b_{\mathrm{irr}}=2^{2g-6}$ and $b_{g-1}=2^{g-3}(2^{g-1}+1)$. This completes the proof of Theorem \ref{oddthet}.

\end{document}